% To: shlhetal@math.huji.ac.il
% Date: Thu, 22 Jan 1998 15:11:24 +0100
% From: spinas@math.ethz.ch
% X-sliced-and-diced-by: 'savemail' 0.1, Aug 30, 1997

\def\forces{\parallel\!\!\! -}

%Das Folgende sollte die "Blackboard bold font" \Bbb definieren, und 
%auch  die "Fraktur" font, \frak:
% Beispiel:    $ {\frak c} = |{\Bbb R}| $

\def\hexnumber#1{\ifcase#1 0\or1\or2\or3\or4\or5\or6\or7\or8\or9\or
        A\or B\or C\or D\or E\or F\fi }

%  The following lines establish the use of the Euler Fraktur font.
\font\teneuf=eufm10
\font\seveneuf=eufm7
\font\fiveeuf=eufm5
\newfam\euffam
\textfont\euffam=\teneuf
\scriptfont\euffam=\seveneuf
\scriptscriptfont\euffam=\fiveeuf
\def\frak{\fam\euffam \teneuf}
%  End definition of Euler Fraktur font.

\font\tenmsx=msam10
\font\sevenmsx=msam7
\font\fivemsx=msam5
\font\tenmsy=msbm10
\font\sevenmsy=msbm7
\font\fivemsy=msbm5
\newfam\msxfam
\newfam\msyfam
\textfont\msxfam=\tenmsx  \scriptfont\msxfam=\sevenmsx
  \scriptscriptfont\msxfam=\fivemsx
\textfont\msyfam=\tenmsy  \scriptfont\msyfam=\sevenmsy
  \scriptscriptfont\msyfam=\fivemsy
\edef\msx{\hexnumber\msxfam}

\mathchardef\upharpoonright="0\msx16

\def\Bbb#1{\tenmsy\fam\msyfam#1}

\def\Smallskip{\vskip1.4truecm}
\def\Bigskip{\vskip2.2truecm}

\def\qed{{\vcenter{\hrule height.4pt \hbox{\vrule width.4pt height5pt
 \kern5pt \vrule width.4pt} \hrule height.4pt}}}
\def\notin{{\in}\kern-5.5pt / \kern1pt}
\def\ok{\vbox{\hrule height 8pt width 8pt depth -7.4pt
    \hbox{\vrule width 0.6pt height 7.4pt \kern 7.4pt \vrule width 0.6pt height 7.4pt}
    \hrule height 0.6pt width 8pt}}
\def\nt{{\leq}\kern-1.5pt \vrule height 6.5pt width.8pt depth-0.5pt \kern 1pt}
\def\sd{{\times}\kern-2pt \vrule height 5pt width.6pt depth0pt \kern1pt}
\def\zp#1{{\hochss Y}\kern-3pt$_{#1}$\kern-1pt}

\def\extend { \hat {\; \; } }

\def\LL{{\Bbb L}}

\def\MM{{\Bbb M}}

\def\sm{{\smallskip}}

\def\la{{\langle}}
\def\ra{{\rangle}}
\def\sub{\subseteq}
\def\alm{\sub ^*}

\def \o {\omega }
\def \fun {{^\omega \omega }}

\def\seq{{^{<\omega }\omega }}

\font\capit=cmcsc10 scaled\magstep0

\font\bolds=cmssdc10 scaled\magstep0

\overfullrule=0pt
\openup1.5\jot

\def \inf {[\o ]^\o }

\def \o {\omega } \def \sub {\subseteq }
\def \fun {{^{\omega }\omega }}
\def \k {\kappa } \def \lam {\lambda }
\def \l {\langle } \def \r {\rangle }
\def \g {\gamma } \def \s {\sigma } \def\t {\tau } \def \c {{\frak c}}
\font\gross=cmbx10 scaled \magstep1
\font\sgross=cmbx10 scaled \magstep2
  
 \def \P {{\cal P}} \def \A {{\cal A}}
\def\Q {Q(D_0,\dots ,D_{n-1})}
\def\G {G(D_0,\dots ,D_{n-1})}

\def\a {\alpha } \def\b {\beta } \def \d {\delta } \def \h {{\frak h}}
\noindent {\sgross The distributivity numbers of finite products of
$\P (\o )$/fin}

\Bigskip 

\noindent Saharon Shelah\footnote{$^1$}{The author is supported by the
Basic Research Foundation of the Israel Academy of Sciences;
publication 531.}

\smallskip

\item{}{Department of Mathematics, Hebrew University, Givat Ram, 91904
Jerusalem, 

ISRAEL}

\smallskip

\noindent Otmar Spinas\footnote{$^2$}{The author is supported by the
Swiss National Science Foundation}

\smallskip

Mathematik, ETH-Zentrum, 8092 Z\"urich, SWITZERLAND

\Bigskip 

{\narrower

{ABSTRACT: Generalizing [ShSp], for every $n<\o $ we construct a
ZFC-model where the distributivity number of r.o.$(\P (\o
)/\hbox{fin})^{n+1}$, $\h (n+1)$, is smaller than the one of r.o.$(\P (\o
)/\hbox{fin})^{n}$. This answers an old problem of Balcar, Pelant and
Simon (see [BaPeSi]). We also show that both of Laver and Miller forcing
collapse the continuum to $\h (n)$ for every $n<\o $, hence by the
first result, consistently they collapse it below $\h (n)$.

}}

\Bigskip

{\gross Introduction}

\Bigskip

For $\lam $ a cardinal let $\h (\lam )$ be the least cardinal $\k $ for
which r.o.$(\P (\o ) $/fin$)^\lam $ is not $\k -$distributive, where by 
$(\P (\o )$/fin$)^\lam $ we mean the (full) $\lam -$product of $\P (\o
)$/fin in the
forcing sense; so $f\in (\P (\o )$/fin$)^\lam $ if and only if $f: \lam
\rightarrow \P (\o )$/fin $\setminus \{ 0\} $, and the ordering is
coordinatewise.

In [ShSp] the consistency of $\h (2)<\h$ with ZFC has been proved,
which provided a (partial) answer to a question of Balcar, Pelant and
Simon in [BaPeSi]. This inequality holds in a model obtained by
forcing with a countable support iteration of Mathias forcing over a
model of GCH. The proof is long and difficult. The following are the
key properties of Mathias forcing (M.f.) which are essential to the
proof (see [ShSp] or below for precise definitions):

\sm

\item{(1)} M.f. factors into a $\s$-closed and a $\s$-centered
forcing.

\item{(2)} M.f. is Suslin-proper which means that, firstly, it is
simply definable, and, secondly, it permits generic conditions over
every countable model of ZF$^-$.

\item{(3)} Every infinite subset of a Mathias real is also a Mathias
real.

\item{(4)} Mathias forcing does not change the cofinality of any
cardinal from above $\h$ to below $\h$.

\item{(5)} Mathias forcing has the pure decision property and it has
the Laver property.

\sm

In this paper we present a forcing $Q^n$, where $0<n<\o$, which is an
$n$-dimensional version of M.f. which satisfies all the analogues of
the five key properties of M.f. In this paper we only prove
these. Once this has been done the proof of [ShSp] can be generalized
in a straightforward way, to prove the following:

\sm

{\bolds Theorem.} {\it Suppose $V\models ZFC +GCH$. If $P$ is a countable
support iteration of $Q^{n}$ of length $\o_2$ and $G$ is $P$-generic
over $V$, then $V[G]\models \h(n+1)=\o _1 \wedge \h(n)= \o_2$.}

\sm

Besides the fact that the consistency of $\h (n+1)<\h (n)$ was an open
problem in [BaPeSi], our motivation for working on it was that in
[GoReShSp] it was shown that both of Laver and Miller forcing collapse the
continuum to $\h $. Moreover, using ideas from [GoJoSp] and [GoReShSp]
it can be proved that these forcings do not collapse $\c $ below $\h
(\o )$. We do not know whether they do collapse it to $\h (\o )$. But
in $\S 2$ we show that they collapse it to $\h (n)$, for every $n<\o
$. Combining this with the first result we conclude that for every
$n<\o $, consistently Laver and Miller forcing collapse $\c $ strictly
below $\h (n)$.

The reader should have a copy of [ShSp] at hand. We
do not repeat all the definitions from [ShSp] here. Notions as Ramsey
ultrafilter, Rudin-Keisler ordering, Suslin-proper are explained there
and references are given. 

\Bigskip 

{\gross 1. The forcing}

\Bigskip

{\bolds Definition 1.1.}  Suppose that $D_0, \dots , D_{n-1}$ are
ultrafilters on $\o $. The game $G(D_0,\dots ,$ $D_{n-1})$ is defined as
follows: In his $m$th move player I chooses $\la A_0,\dots ,A_{n-1}\ra \in
D_0\times \dots \times D_{n-1}$ and player II responds playing $k_m\in
A_{m\hbox{mod}n}$. Finally player II wins if and only if for every
$i<n$, $\{ k_j : j=i\hbox{mod}n \} \in D_i$ holds. 

\sm

{\bolds Lemma 1.2.} {\it Suppose $D_0,\dots ,D_{n-1}$ are Ramsey
ultrafilters which are pairwise not RK-equivalent. Let $\langle
m(l):l<\o \rangle$ be an increasing sequence of integers. There exists
a subsequence $\langle m(l_j): j<\o\rangle $ and sets $Z_i\in D_i$,
$i<n$, such that:

\item{(1)} $l_{j+1}-l_j\geq 2$, for all $j<\o$,}

\item{(2)} $Z_i\sub \bigcup_{j=i\hbox{mod}n}[m(l_j), m(l_{j+1}))$,
{\it for all} $i<n$,

\item{(3)} $Z_i\cap [m(l_j),m(l_{j+1}))$ {\it has precisely one member, for
every $i<n$ and} $j=i$mod$n$.

\sm

{\capit Proof:} For $j<3$, $k<\o$ define:
$$I_{j,k}=\bigcup_{s=(2n-1)(3k+j)}^{(2n-1)(3k+j+1)-1} [m_s, m_{s+1})\,
,$$ $$J_j=\bigcup_{k<\o }I_{j,k}\, .$$ As the $D_i$ are Ramsey
ultrafilters, there exist $X_i\in D_i$ such that for every $i<n$:

\item{(a)} $X_i\sub J_j$ for some $j<3$,

\item{(b)} if $X_i\sub J_j$, then $X_i\cap I_{j,k}$ contains precisely
one member, for every $k<\o$.

Next we want to find $Y_i\in D_i$, $Y_i\sub X_i$, such that for every
distinct $i, i' <n$, $Z_i$ and $Z_{i'}$ do not meet any adjacent
intervals $I_{j,k}$.

Define $h:X_0\rightarrow X_1$ as follows. Suppose $X_0\sub J_j$. For
every $k<\o$, $h$ maps the unique element of $X_0\cap I_{j,k}$ to the
unique element of $X_1$ which belongs to either $I_{j,k}$ or to one of
the two intervals of the form $I_{j',k'}$ which are adjacent to
$I_{j,k}$ (note that these are $I_{2,k-1}, I_{1,k} $ if $j=0$, or
$I_{0,k}, I_{2,k}$ if $j=1$, or $I_{1,k}, I_{0,k+1}$ if $j=2$). As $h$
does not witness that $D_0, D_1$ are RK-equivalent, there exist
$X_i'\in D_i$, $X_i'\sub X_i$ ($i<2$) such that $h[X_0']\cap
X_1'=\emptyset$. Note that if $n=2$, we can let $Y_i =X_i'$. Otherwise
we repeat this procedure, starting from $X_0'$ and $X_2$, and get
$X_0''$ and $X_2'$. We repeat it again, starting from $X_1'$ and
$X_2'$, and get $X_1''$ and $X_2''$. If $n=3$ we are done. Otherwise
we continue similarly. After finitely many steps we obtain $Y_i$ as
desired. 

By definition of $I_{j,k}$ it is now easy to add more elements to each
$Y_i$ in order to get $Z_i$ as in the Lemma. The ``worst'' case is
that some $Y_i$ contains integers $s<t$ such that $(s,t)\cap
Y_{u}=\emptyset $ for all $u <n$. By construction there is some
$I_{j,k}\sub (s,t)$. For every $u<n-1$ pick $$x_u\in
[m((2n-1)(3k+j)+2u+1), m((2n-1)(3k+j)+2u+2))$$ and add $x_u$ to
$Y_{i+u+1\hbox{mod}n}$. The other cases are similar. $\qed$

\sm

{\bolds Corollary 1.3.}  {\it Suppose $D_0,\dots ,D_{n-1}$ are Ramsey
ultrafilters which are pairwise not RK-equivalent. Then in the game
$G(D_0,\dots ,D_{n-1})$ player I does not have a winning strategy.}

\sm

{\capit Proof:} Suppose $\s $ is a strategy for player I. For every
$m<\o , i<n$ let $\A ^m_i\sub D_i$ be the set of all $i$th coordinates
of moves of player I in an initial segment of length at most $2m+1$ of
a play in which player I follows $\s $ and player II plays only
members of $m$.

As the $D_i$ are $p-$points and each $\A ^m_i$ is finite, there exist
$X_i\in D_i$ such that $\forall m\forall i<n\forall A\in \A
^m_i(X_i\alm A)$. Moreover we may clearly find a strictly increasing
sequence $\la m(l):l<\o \ra $ such that $m(0)=0$ and for all $l<\o $:

$$\forall i<n\forall A\in \A ^{m(l)}_i(X_i \sub A\cup m({l+1}) \wedge
X_i \cap [m(l),m({l+1})) \ne \emptyset ).$$

Applying Lemma 1.2, we obtain a subsequence $\la m(l_j):j<\o\ra$ and
sets $Z_i\in D_i$. 

Now let in his $j$th move player II play $k_j$, where $k_j$ is the
unique member of $[m({l_j}), m({l_{j+1}}))\cap X_{j\hbox{mod}n}\cap
Z_{j\hbox{mod}n}$ if it exists, or otherwise is any member of
$[m({l_j}), m({l_{j+1}})) \cap X_{j\hbox{mod}n}$ (note that this
intersection is nonempty by definition of $m(l_{j+1})$. Then this play is
consistent with $\s $, moreover $X_i\cap Z_i \sub \{ k_j :j=i\hbox{mod}n\} $
for every $i<n$, and hence it is won by player II. Consequently $\s $
could not have been a winning strategy for player I. $\qed $
 
\sm

{\bolds Definition 1.3.}  Let $n<\o $ be fixed. The forcing $Q$ (really $Q^n$) is defined as follows: Its
members are $(w, \bar A)\in [\o ]^{<\o }\times [\o ]^\o $. If
$\la k_j :j<\o \ra $ is the increasing enumeration of $\bar
A$ we let $\bar A_i =\{ k_j : j=i\hbox{mod}n \} $ for $i<n$, and if
$\la l_j:j<m \ra $ is the increasing enumeration of $w$ then
let $w_i = \{ l_j:j=i\hbox{mod}n\} $, for $i<n$.

Let $(w,\bar A)\leq (v,\bar B)$ if and only if $w\cap (\max (v)+1)=v$,
$w_i\setminus v_i \sub \bar B_i$ and $\bar A_i \sub \bar B_i$, for
every $i<n$.

If $p\in Q$, then $w^p, w_i^p, \bar A^p, \bar A^p_i$
have the obvious meaning. We write $p\leq ^0q$ and say ``$p$ is a pure
extension of $q$'' if $p\leq q$ and $w^p=w^q$.

If $D_0, \dots ,D_{n-1}$ are ultrafilters on
$\o $, let $\Q $ denote the subordering of $Q$ containing only those
$(w, \bar A)\in Q$ with the property $\bar A_i\in D_i$, for every $i<n$.

\sm

{\bolds Lemma 1.4.}  {\it The forcing $Q$ is equivalent to $(\P (\o
)\hbox{/fin})^n \ast Q(\dot G_0, \dots \dot G_{n-1})$, where $(\dot
G_0,\dots ,\dot G_{n-1})$ is the canonical name for the generic
object added by $(\P (\o )\hbox{/fin})^n$, which consists of $n$
pairwise not RK-equivalent Ramsey ultrafilters.}

\sm

{\capit Proof:}  Clearly $(\P (\o )\hbox{/fin})^n$ is $\s
-$closed and hence does not add reals. Moreover, members $\la
x_0,\dots ,x_{n-1}\ra \in (\P (\o )\hbox{/fin})^n$ with the property
that if $\bar A=\bigcup \{ x_i:i<n\} $, then $x_i=\bar A_i$ for every $i<n$
are dense. Hence the map $(w,\bar A) \mapsto (\la \bar A_0,\dots
,\bar A_{n-1}\ra , (w,\bar A))$ is a dense embedding of the respective
forcings. 

That $\dot G_0,\dots ,\dot G_{n-1}$ are ($(\P (\o
)\hbox{/fin})^n-$forced to be) pairwise not
RK-equivalent Ramsey ultrafilters follows by an easy genericity
argument and again the fact that no new reals are added. $\qed $ 

\sm

{\bolds Notation.}  We will usually abbreviate the decomposition of
$Q$ from Lemma 1.4. by writing $Q=Q'\ast Q''$. So members of
$Q'$ are $\bar A,\bar B\in \inf $ ordered by $\bar A_i\sub \bar B_i$ for all
$i<n$; $Q''$ is $Q(\dot G_0,\dots ,\dot G_{n-1})$. If $G$ is a
$Q-$generic filter, by $G'\ast \dot G''$ we denote its decomposition
according to $Q=Q'\ast \dot Q''$, and we write $G'= (G'_0,\dots ,G'_{n-1})$.

\sm

{\bolds Definition 1.5.}  Let $I\sub \Q $ be open dense. 
%and $p\in \Q $. 
We define a rank function rk$_{I}$ on $[\o ]^{<\o }$
%the set of all $w\in [\o
%]^{<\o }$ with $(w,\bar A^p)\leq p$ 
as follows. Let rk$_{I}(w)=0$ if
and only if $(w,\bar A)\in I$ for some $\bar A$. Let rk$_{I}(w)=\a $ if
and only if $\a $ is minimal such that there exists $A\in
D_{|w|\hbox{mod}n} $ with the property that for every $k\in A$,
rk$_{I}(w\cup \{ k\} )=\b $ for some $\b <\a $. Let rk$_I(w)=\infty $
if for no ordinal $\a $, rk$_I(w)=\a $.

\sm

{\bolds Lemma 1.6.}  {\it If $D_0,\dots ,D_{n-1}$ are Ramsey
ultrafilters which are pairwise not RK-equivalent and $I\sub \Q $ is
open dense, then for every $w\in [\o ]^{<\o }$, rk$_I(w)\ne \infty $.}

\sm

{\capit Proof:} Suppose we had rk$_I(w)=\infty $ for some $w$. We
define a strategy $\s $ for player I in $\G $ as follows: $\s
(\emptyset )=\la A_0,\dots ,A_{n-1}\ra \in D_0 \times \dots \times
D_{n-1}$ such that for every $k\in A_{|w|\hbox{mod}n}$, rk$_I(w \cup
\{ k\} )=\infty $. This choice is possible by assumption and as the
$D_i$ are ultrafilters. In general, suppose that $\s $ has been
defined for plays of length $2m$ such that whenever $k_0,\dots
,k_{m-1}$ are moves of player II which are consistent with $\s $, then
$k_0< k_1\dots < k_{m-1}$ and for every $\{ k_{i_0} < \dots
<k_{i_{l-1}}\} \sub \{ k_0,\dots ,k_{m-1}\} $ with $i_j=j\hbox{mod}n,
j<l,$ we have rk$_I(w\cup \{ k_{i_0},\dots ,k_{i_{l-1}}\})=\infty
$. Let $S$ be the set of all $\{ k_{i_0} < \dots <k_{i_{l-1}}\} \sub
\{ k_0,\dots ,k_{m-1}\} $ with $i_j=j\hbox{mod}n, j<l,$ and
$l=m$mod$n$. As $D_{|w|+m{\hbox{mod}}n}$ is an ultrafilter, by
induction hypothesis we have that, letting $$A_{|w|+m{\hbox{mod}}n}
=\{ k>k_{m-1}: \forall s\in S\quad \hbox{rk}_I(w\cup s\cup \{
k\})=\infty \},$$ $A_{|w|+m{\hbox{mod}}n}\in
D_{|w|+m{\hbox{mod}}n}$. For $i\ne |w|+m$mod$n$, choose $A_i\in D_i$
arbitrarily, and define $$\s \la k_0,\dots k_{m-1}\ra =\la A_0,\dots
,A_{n-1}\ra. $$ 

Since by Lemma 1.2. $\s $ is not a winning strategy for player I,
there exist $k_0<\dots <k_m< \dots $ which are moves of player II
consistent with $\s $, such that, letting $\bar A=\{ k_m:m<\o \} $,
we have $(w,\bar A)\in \Q $. By construction we have that for every
$(v,\bar B)\leq (w,\bar A)$, rk$_I(v)=\infty $.
This contradicts the assumption that $I$ is dense.

\sm

{\bolds Definition 1.7.}  Let $p\in Q$. A set of the form $w^p\cup \{
k_{|w|}<k_{|w|+1}<\dots \}\in \inf $ is called a {\it branch} of $p$
if and only if $\max (w^p)< k_{|w|}$ and $\{ k_j:j=i\hbox{mod}n\} \sub
\bar A^p_i$, for every $i<n$. A set $F\sub [\o ]^{<\o }$ is called a
{\it front} in $p$ if for every $w\in F$, $(w, \bar A^p)\leq p$ and
for every branch $B$ of $p$, $B\cap m\in F$ for some $m<\o $.

\sm

{\bolds Lemma 1.8.}  {\it Suppose $D_0,\dots , D_{n-1}$ are pairwise
not RK-equivalent Ramsey ultrafilters. Suppose $p\in \Q $ and $\la
I_m:m<\o \ra $ is a family of open dense sets in $\Q $. There exists
$q\in \Q $, $q\leq ^0p$, such that for every $m$, $\{ w\in [\o ]^{<\o
}: (w, \bar A^q )\in I_m \wedge (w,\bar A^q)\leq q\} $ is a front in $q$.}

\sm

{\capit Proof:} First we prove it in the case $I_m=I$ for all $m<\o $,
by induction on rk$_I(w^p)$. We define a strategy $\s $ for player I
in $\G $ as follows. Generally we require that $$\s\la k_0,\dots ,
k_r\ra _i \sub \s \la k_0,\dots ,k_s\ra _i$$ for every $s<r$ and
$i<n$, where $\s\la k_0,\dots , k_r\ra _i$ is the $i$th coordinate of
$\s\la k_0,\dots , k_r\ra$. We also require that $\s$ ensures that the
moves of II are increasing (see the proof of 1.7). Define $\s
(\emptyset)=\la A_0,\dots ,A_{n-1}\ra $ such that for every $k\in
A_{|w^p|\hbox{mod}n }$, rk$_I(w^p\cup \{ k\})< $ rk$_I(w^p)$. 

Suppose now that $\s$ has been defined for plays of length $2m$, and
let $\la k_0,\dots ,k_{m-1}$ be moves of II, consistent with $\s$. The
interesting case is that $m-1 =0$mod$n$. Let us assume this first. By
definition of $\s (\emptyset )$ and the general requirement on $\s$ we
conclude rk$_I(w^p\cup \{ k_{m-1}\}) <$ rk$_I(w^p)$. By induction
hypothesis there exists $\la A_0,\dots ,A_{n-1}\ra\in D_0\times \dots
\times D_{n-1}$ such that, letting $\bar A=\bigcup_{i<n}A_i$, we have
$(w^p, \bar A)\leq p$ and $$\{ v\in [\o]^{<\o}: (v,\bar A)\in I \wedge
(v,\bar A)\leq (w^p\cup \{ k_{m-1}\}, \bar A)\} $$ is a front in
$(w^p\cup \{ k_{m-1}\}, \bar A) $. We shrink $\bar A$ such that,
letting $$\s\la k_0,\dots ,k_{m-1}\ra =\la A_0,\dots ,A_{n-1}\ra ,$$
the general requirements on $\s$ above are satisfied. 

In the case that $m-1 \ne 0$mod$n$, define $\s\la k_0,\dots
,k_{m-1}\ra$ arbitrarily, but consistent with the rules and the
general reqirements above.

Let $\bar A =\{ k_i :i<\o\}$ be moves of player II witnessing that
$\s$ is not a winning strategy. Let $q=(w^p, \bar A)$. Let $B= w^p
\cup \{ l_{|w^p|}< l_{|w^p|+1} < \dots \}$ be a branch of $q$. Hence
$l_{|w^p|}= k_j$ for some $j=0$mod$n$. Then $w^p\cup \{ k_j\} \cup \{
l_{|w^p|+1} , l_{|w^p|+2}, \dots \}$ is a branch of $(w^p \cup \{
k_j\} , \s\la k_0, \dots , k_j\ra )$. By definition of $\s$ there
exists $m$ such that $(B\cap m, \s\la k_0, \dots , k_j\ra)\in I$. As
$(B\cap m, \bar A)  \leq (B\cap m, \s\la k_0, \dots , k_j\ra)$ and $I$
is open we are done. 

\smallskip

For the general case where we have infinitely many $I_m$, we make a
diagonalization, using the first part of the present proof. Define a
strategy $\s $ for player I satisfying the same general requirements
as in the first part as follows. Let $\s (\emptyset )= \la A_0,\dots
,A_{n-1}\ra $ such that, letting $\bar A=\bigcup \{ A_i:i<n\} $,
$(w^p,\bar A)\leq ^0p$ and it satisfies the conclusion of the Lemma
for $I_0$. In general, let $\s \la k_0,\dots , k_{m-1}\ra =\la
A_0,\dots ,A_{n-1}\ra $ such that, letting $\bar A=\bigcup \{
A_i:i<n\} $, for every $v\sub \{ k_i:i<m\} $ and $j\leq m$, $(w^p\cup
v,\bar A)\leq ^0 (w^p\cup v,\bar A^p) $ and it satisfies the
conclusion of the Lemma for $I_j$ (In fact we don't have to consider
all such $v$ here, but it does not hurt doing it). If then $\bar A=\{
k_i:i<\o \} $ are moves of player II witnessing that $\s $ is not a
winning strategy for I, similarly as in the first part it can be
verified that $q=(w^p,\bar A)$ is as desired. $\qed $

\sm

{\bolds Corollary 1.9.}  {\it Let $D_0,\dots ,D_{n-1}$ be pairwise not
RK-equivalent Ramsey ultrafilters. Suppose $\bar A\in \inf $ is such
that for every $i<n$ and $X\in D_i $, $\bar A_i \alm X$. Then $\bar A$ is
$\Q -$generic over $V$.}

\sm

{\capit Proof:}  Let $I\sub \Q $ be open dense. Let $w\in [\o ]^{<\o
}$. It is easy to see that the set 

$$I_w=\{ (v,\bar B)\in \Q :  (w\cup [v \setminus \min \{ k \in
v_{|w|\hbox{mod}n }: k>\max (w)\} ],\bar B )\in I \} $$

\noindent is open dense. If we apply Lemma 1.8. to $p=(\emptyset ,\o
,\dots ,\o )$ and the countably many open dense sets $I_w$ where $w\in
[\o ]^{<\o }$, we obtain $q=(\emptyset ,\bar B)$. Let $\la a_i:i<\o
\ra $ be the increasing enumeration of $\bar A$. Choose $m$ large
enough so that for each $i<n$, $\bar A_i \setminus \{ a_j: j<mn\} \sub
\bar B_i $. Let $w=\{ a_j:j<mn \} $. By construction, there exists
$v\sub \bar A\cap \bar B\setminus (a_{mn -1}+1)$ such that $(v,\bar
B)\in I_w$ and $w\cup v=\bar A \cap k$, for some $k<\o $.  Hence
$(w\cup v, \bar B)\in I$, and so the filter on $\Q $ determined by
$\bar A$ intersects $I$. As $I$ was arbitrary, we are done. $\qed $

\sm

An immediate consequence of Lemma 1.4. and Corollary 1.9. is the
following. 

\sm

{\bolds Corollary 1.10.}  {\it Suppose $\bar A\in \inf $ is $Q-$generic over
$V$, and $\bar B\in \inf $ is such that $\bar B_i\sub \bar A_i$
for every $i<n$.
Then $\bar B$ is $Q-$generic over $V$ as well.} 

\sm

Remember that a forcing is called Suslin, if its underlying set is an
analytic set of reals and its order and incompatibility relations are
analytic subsets of the plane. A forcing $P$ is called Suslin-proper
if it is Suslin and for every countable transitive model $(N,\in )$ of
ZF$^-$ which contains the real coding $P$ and for every $p\in P\cap
N$, there exists a $(N,P)-$generic condition extending $p$. See [JuSh]
for the theory of Suslin-proper forcing and [ShSp] for its
properties which are relevant here.

\sm

{\bolds Corollary 1.11.}  {\it The forcing $Q$ is Suslin-proper.}

\sm

{\capit Proof:}  It is trivial to note that $Q$ is Suslin, without
parameter in its definition. Let $(N,\in)$ be a countable model of
ZFC$^-$, and let $p\in Q\cap N$. Without loss of generality,
$|w^p|=0\hbox{mod} n$. Let $\bar A\in \inf \cap V$ be $Q-$generic over
$N$ such that $p$ belongs to its generic filter. Hence $w^p_i \sub
\bar A_i \sub w^p_i \cup (\bar A^p_i \setminus (\max (w^p)+1))$ for
all $i<n$. But if
$q=(w^p,\bar A)$, then clearly $q\leq ^0p$ and $q$ is $(N,Q)-$generic,
as every $\bar B\in \inf $ which is $Q-$generic over $V$ and contains
$q$ in its generic filter is a subset of $\bar A$ and hence $Q\cap
N-$generic over $N$ by Corollary 1.10. applied in $N$. $\qed $

\sm

The following is an immediate consequence of Corollary 1.11.

\sm

{\bolds Corollary 1.12.}  {\it If $p\in Q$ and $\la \t _n:n<\o \ra $
are $Q-$names for members of $V$, there exist $q\in Q$, $q\leq ^0p$
and $\la X_n:n<\o \ra $ such that $X_n\in V\cap [V]^\o $ and $q\forces
_Q\; \forall n (\t _n\in X_n)$.}

\sm

{\bolds Corollary 1.13.}  {\it Forcing with $Q$ does not change the
cofinality of any cardinal $\lambda $ with cf$(\lambda )\geq \h (n)$
to a cardinal below $\lambda $.}

\sm

{\capit Proof:}  Suppose there were a cardinal $\kappa < \h (n)$ and a
$Q-$name $\dot f$ for a cofinal function from $\kappa $ to $\lambda $.
Working in $V$ and using Corollary 1.12., for every $\a <\kappa $ we
may construct a maximal antichain $\la p^\a _\b :\b < {\frak c} \ra $
in $Q$ and $\la X^\a _\b :\b <{\frak c}\ra $ such that for all $\b <\c
$, $w^{p^\a _\b }=\emptyset $, $X^\a _\b \in [V]^\o \cap V$ and $p^\a
_\b \forces _Q\; \dot f (\a )\in X^\a _\b $.

Then clearly $\A _\a = \la \la    \bar A_i^{p^\a _\b }:i<n\ra :\b <\c
\ra $ is a maximal antichain in $(\P (\o )$/fin$)^n$. By $\kappa <\h
(n)$, $\la \A _\a :\a <\kappa \ra $ has a refinement, say $\A $.
Choose $\la \bar A_i:i<n\ra \in \A $. Let $\bar A=\bigcup \{ \bar A_i:i<n \} $.
We may assume that the $\bar A_i$ also have the meaning from Definition
1.3. with respect to $\bar A$. For each $\a <\kappa $ there exists $\b
(\a )$ such that $\la \bar A_i:i<n\ra \leq _{(\P (\o )/fin)^n } \la
\bar A_i^{p^\a _{\b (\a )}}: i<n\ra $. Then clearly 

$$(\emptyset ,\bar A)\forces _Q\; \hbox{range}(\dot f)\sub \bigcup \{
X^\a _{\b (\a )}:\a <\kappa \} \, .$$

But as cf$(\lambda )\geq \h (n)$ and $\kappa <\h (n)$, we have a
contradiction. $\qed $

\sm

{\bolds Lemma 1.14.}  {\it Suppose $D_0,\dots ,D_{n-1}$ are pairwise
not RK-equivalent Ramsey ultrafilters. Then $\Q $ has the pure
decision property (for finite disjunctions), i.e. given a $\Q -$name
$\t $ for a member of $\{ 0,1\} $ and $p\in \Q $, there exist $q\in
\Q $ and $i\in \{ 0,1\} $ such that $q\leq ^0p$ and $q\forces _{\Q }\;
\t =i $.}

\sm

{\capit Proof:} The set $I=\{ r\in \Q : r$ decides $ \t \} $ is open
dense. By a similar induction on rk$_I$ as in the proof of Lemma 1.7.
we may find $q\in \Q $, $q\leq ^0p$, such that for every $q'\leq q$,
if $q'$ decides $ \t $ then $(w^{q'}, \bar A^q)$ decides $\t $. Now again
by incuction on rk$_I$ we may assume that for every $k\in \bar
A^q_{|w^q|\hbox{mod}n} $, $(w^q\cup \{ k\} ,\bar A^q)$ satisfies the
conclusion of the Lemma, and hence by the construction of $q$,
$(w^q\cup \{ k\} ,\bar A^q)$ decides $ \t $. But then clearly a pure
extension of $q$ decides $\t $, and hence $q$ does. $\qed $

\sm

{\bolds Lemma 1.15.} {\it Lemma 1.14 holds if $\Q $ is replaced by
$Q$.} 

\sm

{\capit Proof:} Suppose $p\in Q$, $\t $ is a $Q-$name and $p\forces
_Q\; \t \in \{ 0,1\} $. As $\bar A^p\forces _{Q'} \; ``p\in Q(\dot
G_0,\dots ,\dot G_{n-1})$'', by Lemma 1.14 there exists a $Q'-$name
$\dot {\bar A}$ such that 

$$\bar A^p\forces _{Q'}\; ``(w^p, \dot {\bar A})\in Q''
\wedge (w^p,\dot{\bar A})\leq p \wedge (w^p,\dot {\bar
A})\hbox{ decides } \t \hbox{''}\, .$$

As $Q'$ does not add reals there exist $\bar A_1,\bar A_2\in \inf \cap
V$ such that $\bar A_1\sub \bar A^p$ and $\bar A_1 \forces _{Q'} \;
\dot {\bar A} =\bar A_2$. Letting $\bar B =\bar A_1\cap \bar A_2$ we
conclude $(w^p,\bar B)\in Q$, $(w^p,\bar B)\leq ^0p$ and $(w^p,\bar
B)\hbox{ decides } \t $. $\qed $

\sm

%{\bolds Lemma 1.14.}  {\it Suppose $D_0,\dots ,D_{n-1}$ are pairwise
%not RK-equivalent Ramsey ultrafilters, $p\in \Q $ and $\la \t _n:n<\o
%\ra $ are $\Q -$names for members of $V$. Then there exists $q\in \Q
%$, $q\leq ^0$, such that for every $m<\o $ and $q'\in \Q $ with
%$q'\leq q$ and $w^{q'}\sub m$, for every $i\leq m$, $q'||_{\Q } \; \t
%_i$ if and only if $(w^{q'}, \bar A^q\setminus \min (\bar
%A^q_0\setminus m))||_{\Q }\; \t _i $.}
%{\capit Proof:} Similarly as in several previous Lemmas define an
%appropriate strategy $\s $ for player I in $\G $ such that if $\bar A= \{
%k_0<\dots <k_i<\dots \} $ is the set of moves of player II witnessing
%that $\s $ is not a winning strategy, then $q=(w^p,\bar A)$ is as
%desired. $\qed $
%{\bolds Corollary 1.15.}  {\it The assertion of Lemma 1.14. holds if
%$\Q $ is replaced by $Q$.}
%{\capit Proof:}  The proof is similar as the one of Lemma 1.13. $\qed
%$ 

\sm

The rest of this section is devoted to the proof that if the forcing
$Q$ is iterated with countable supports, then in the resulting model
cov$({\cal M})=\o _1$, where ${\cal M}$ is the ideal of meagre subsets
of the real line, and cov$({\cal M})$ is the least number of meagre sets
needed to cover the real line. Hence for every $n<\o $, 
we obtain the consistency of cov$({\cal M})< \h (n)$. 

\sm

{\bolds Definition 1.16.}  A forcing $P$ is said to have the {\it
Laver property} if for every $P-$name $\dot f$ for a member of $\fun
$, $g\in \fun \cap V$ and $p\in P$, if 

$$p\forces _P\; \forall n<\o (\dot f(n)< g(n))\, ,$$

\noindent then there exist $H:\o \rightarrow [\o ]^{<\o }$ and $q\in
P$ such that $H\in V$, $\forall n<\o (|H(n)|\leq 2^n)$, $q\leq p$ and

$$q\forces _P \; \forall n<\o (\dot f(n)\in H(n))\, .$$

\sm

It is not difficult to see that a forcing with the Laver property does
not add Cohen reals. Moreover, by [Shb, 2.12., p.207] the Laver
property is preserved by a countable support iteration of proper
forcings. See also [Go, 6.33., p.349] for a more accessible proof. 

\sm

{\bolds Lemma 1.17.}  {\it The Forcing $Q$ has the Laver property.}

\sm

Suppose $\dot f$ is a $Q-$name for a member of $\fun $ and $g\in \fun
\cap V$ such that $p\forces _Q \, \forall n<\o (\dot f(n)<g(n))$. We
shall define $q\leq ^0p$ and $\l H(i):i<\o \r $ such that $|H(i)|\leq
2^i$ and $q\forces _Q\, \forall i(\dot f(i)\in H(i))$. We may assume
$|w^p |=0$mod$n$ and $\min (\bar A^p)>\max (w^p)$. 

By Lemma 1.14 choose $q_0\leq ^0p$ and $K^0$ such that $q_0\forces
_Q\, \dot f(0)=K^0$, and let $H(0)=\{ K^0\} $.

Suppose $q_i \leq ^0p$, $\l H(j):j\leq i\r $ have been constructed
and let $a^i$ be the set of the first $i+1$ members of $\bar A^{q_i}$.
Let $\l v^k:k<k^*\r $ list all subsets $v$ of $a^i$ such that $v_l\sub
(a^i)_l $, for every $l<n$ (see Definition 1.3.). Then clearly
$k^*\leq 2^{i+1}$. By Lemma 1.14 we may shrink $\bar A^{q_i}$ $k^*$
times so to obtain $\bar A$ and $\l K^{i+1}_k:k<k^* \r $ such that for
every $k<k^*$, $( w^{q_i}\cup v^k, \bar A)\forces _Q\, \dot
f(i+1)=K^{i+1}_k$. Without loss of generality, $\min (\bar A)>\max
(a^i)$. Let $q_{i+1}$ be defined by $w^{q_{i+1}}=w^p$ and $\bar
A^{q_{i+1}}=a^i\cup \bar A'$, where $\bar A'$ is $\bar A$ without its
first $(i+1)$mod$n$ members. Let $H(i+1)=\{ K^{i+1}_k:k<k^* \} $. Then
$q^{i+1}\forces _Q\, \dot f (i+1)\in H(i+1)$. Finally let $q$ be
defined by $w^q=w^p$ and $\bar A^q=\bigcup \{ a^i:i<\o \} $. Then $q$
and $\l H(i):i<\o \r $ is as desired. \hfill $\qed $

\sm

As explained above, from Lemma 1.17 and Shelah's preservation theorem
it follows that if $P$ is a countable support iteration of $Q$ and $G$
is $P-$generic over $V$, then in $V[G]$ no real is Cohen over $V$;
equivalently, the meagre sets in $V$ cover all the reals of $V[G]$.
Now starting with $V$ satisfying CH we obtain the following theorem.

\sm

{\bolds Theorem 1.18.}  {\it For every $n<\o $, the inequality}
cov$({\cal M})< \h (n) $ {\it is consistent with ZFC.}

\Bigskip

{\gross 2. Both of Laver and Miller forcing collapse the continuum below each
$\h (n)$}

\Bigskip

{\bolds Definition 2.1.}  Let $p\sub \seq $ be a tree. For any $\eta
\in p$ let succ$_\eta (p)=\{ n<\o : \eta \extend \la n \ra \in p\} $. We
say that $p$ has a stem and denote it
stem$(p)$, if there is $\eta \in p$ such that
$|$succ$_\eta (p)|\geq 2$ and for every $\nu \subset \eta $,
$|$succ$_\nu (p)|=1$. Clearly stem$(p)$ is uniquely determined, if it
exists. If $p$ has a stem, by $p^-$ we denote the set $\{ \eta \in p:
\hbox{stem}(p) \sub \eta \} $. We say that $p$ is a Laver tree if $p$
has a stem and for every $\eta \in p^-$, succ$_\eta (p)$ is infinite.
We say that $p$ is superperfect if for every $\eta \in p$ there exists
$\nu \in p$ with $\eta \subseteq \nu $ and $|$succ$_\nu (p)|=\o $. By
$\LL $ we denote the set of all Laver trees, ordered by reverse
inclusion. By $\MM $ we denote the set of all superperfect trees,
ordered by reverse inclusion. $\LL ,\MM $ is usually called Laver,
Miller forcing, respectively. 

\sm

{\bolds Theorem 2.2.}  {\it Suppose that $G$ is $\LL -$generic or $\MM
-$generic over $V$. Then in $V[G]$, $|\c ^V|= |\h (n)|^V$.}

\sm

{\capit Proof:} Completely similarly as in [BaPeSi] for the case
$n=1$, a base tree
$T$ for $(\P (\o )$/fin$)^n$ of height $\h (n)$ can be constructed.
I.e. 

\sm

\item{(1)}{$T\sub (\P (\o )$/fin$)^n$ is dense;}

\item{(2)}{$(T,\supseteq ^*)$ is a tree of height $\h (n)$;}

\item{(3)}{each level $T_\a $, $\a <\h (n)$, is a maximal antichain in
$(\P (\o )$/fin$)^n$;}

\item{(4)}{every member of $T$ has $2^\o $ immediate successors.}

\sm   

It follows easily that, firstly, every chain in $T$ of length of
countable cofinality has an upper bound, and secondly, every member of
$T$ has an extension in $T_\a $ for arbitrarily large $\a <\h (n)$. 

Using $T$, we will define a $\LL -$name for a map from $\h (n)$ onto
$\c $. For $p\in \LL $ and $\{ \eta _0,\dots ,\eta _{n-1}\} \in
[p^-]^n$, let $\bar A^p_{\{ \eta _i:i<n\} }=\la \hbox{succ}_{\eta
_i}(p):i<n \ra $.

By induction on $\a <\c $ we will construct $(p_\a ,\d _\a ,\gamma _\a
)\in \LL \times \h (n)\times \c $ such that the following clauses
hold:

\sm

\item{(5)}{if $\{ \eta _0,\dots ,\eta _{n-1 }\} \in [p_\a ]^n$, then
$\bar A^{p_\a }_{\{ \eta _i:i<\o \}}\in T_{\d _\a }$;}

\item{(6)}{if $\beta <\a $, $\d _\b =\d _\a $, $\{ \eta _0,\dots ,\eta
_{n-1}\} \in [p_\a ^-]^n\cap [p_\b ^-]^n$, then $\bar A^{p_\a }_{\{
\eta _i:i<n\} }$, $\bar A^{p_\b }_{\{ \eta _i:i<n\} }$ are
incompatible in $(\P (\o )$/fin$)^n$;}

\item{(7)}{if $p\in \LL$, $\g <\c $, then for some $\a <\c $, every
extension of $p_\a $ is compatible with $p$ and $\g _\a =\g $.}

\sm

At stage $\a $, by a suitable bookkeeping we are given $\g <\c $,
$p\in \LL $, and have to find $\d _\a , p_\a $ such that (5), (6), (7)
hold. For $\eta \in p^-$ let $B_\eta =$ succ$_\eta (p)$; for $\eta \in
\seq \setminus p^-$, $B_\eta =\o $. Let $\la \{ \eta ^i_0,\dots ,\eta
^i_{n-1} \} :i<\o \ra $ list $[\seq ]^n$ such that every member is
listed $\aleph _0$ times.

Inductively we define $\la \xi _i:i<\o \ra $ and $\la B^\rho _\eta :\eta
\in \seq ,\rho \in {^{<\o }2}\ra $ such that 

\sm

\item{(8)}{$B^\rho _\eta \in \inf $ and $\la \xi _i:i<\o \ra $ is a
strictly increasing sequence of ordinals below $\h (n)$;}

\item{(9)}{$B^\emptyset _\eta =B_\eta $;}

\item{(10)}{for every $i<\o $, the map $\rho \mapsto \la B^\rho _{\eta
^i_0},\dots ,B^\rho _{\eta ^i_{n-1}} \ra $ is one-to-one from
$^{i+1}2$ into $T_{\xi _i}$;}

\item{(11)}{for every $i<k$, for every $\rho \in {^{k+1}2}$, $B^\rho
_\eta \sub ^* B^{\rho \upharpoonright i+1}_\eta \alm B^\emptyset _\eta
$;}

\sm

Suppose that at stage $i$ of the construction, $\la \xi _j:j<i\ra $
and $\la B^\rho _\eta : \eta \in \{ \eta ^j_0,\dots ,\eta ^j_{n-1}:j<i
\} , \rho \in {^{\leq i}2} \ra $ have been constructed.

For $\eta \in \{ \eta ^i_0,\dots , \eta ^i_{n-1 } \} $ and $\rho \in
{^{\leq i}2}$, if $B^\rho _\eta $ is not yet defined, there is no
problem to choose it such that (8) and (11) hold. Next by the
properties of $T$ it is easy to find $\xi _i$ and $B^\rho _\eta $, for
every $\rho \in {^{i+1}2}$ and $\eta \in \{ \eta ^i_0, \dots ,\eta
^i_{n-1}\} $ such that $(8), (9), (10), (11)$ hold up to $i$.

By the remark following the properties of $T$, letting $\d _\a =\sup
\{ \xi _i:i<\o \} $, for every $\eta \in \seq $ and $\rho \in {^\o
2}$, there exists $B^\rho _\eta \in \inf $ such that

\sm

\item{(12)}{for all $i<\o $, $B^\rho _\eta \alm B^{\rho
\upharpoonright i}_\eta $;}

\item{(13)}{for all $\{ \eta _0,\dots ,\eta _{n-1}\} \in [\seq ]^n$,
$\la B^\rho _{\eta _0},\dots ,B^\rho _{\eta _{n-1}}\ra \in T_{\d _\a
}$.}

\sm

For $\rho \in {^\o 2}$ let $p^\rho \in \LL $ be defined by 

\sm

stem$(p^\rho )=$ stem$(p_\a )$

$\forall \eta \in (p^\rho )^- ($succ$_\eta (p^\rho )=B^\rho _\eta )$.

\sm

It is easy to see that every extension of $p^\rho $ is compatible with
$p_\a $. Moreover, if $\{ \eta _0,\dots ,\eta _{n-1}\} \in [(p^\rho
)^-]$, then $\bar A^{p^\rho }_{\{ \eta _i:i<n\} }\in T_{\d _\a }$ by
construction. Hence we have to find $\rho \in {^\o 2}$ such that,
letting $p_\a =p^\rho $, (6) holds. Note that for every $\{ \eta
_0,\dots ,\eta _{n-1}\} \in [\seq ]^n$ and $\b <\a $ with $\d _\b =\d
_\a $ and $\{ \eta _0,\dots ,\eta _{n-1 }\} \in [p_\b ^-]^n$ there
exists at most one $\rho \in {^\o 2}$ such that $\{ \eta _0,\dots
,\eta _{n-1 }\} \in [(p^\rho )^-]^n$ and $\bar A^{p^\rho }_{\{ \eta
_i:i<n \} }$, $\bar A^{p_\b }_{\{ \eta _i:i<n\} }$ are compatible in
$(\P (\o )$/fin$)^n$. In fact, by construction and as $T_{\d _\a }$ is
an antichain, either $\bar A^{p^\rho }_{\{ \eta _i:i<n \} } =\bar
A^{p_\b }_{\{ \eta _i:i<n\} }$ or they are incompatible; and moreover
for $\rho \ne \s $, $\bar A^{p^\rho }_{\{ \eta _i:i<n \} }$, $\bar
A^{p^\s } _{\{ \eta _i:i<n\} }$ are incompatible. Hence, as $\aleph _0
\cdot |\a |<\c $ we may certainly find $\rho $ such that, letting $p_\a
=p^\rho $ and $\g _\a =\gamma $, (5), (6), (7) hold.

But now it is easy to define an $\LL -$name $\dot f$ for a function
>from $\h (n)$ to $\c $ such that for every $\a <\c $, $p_\a \forces
_\LL \, \dot f (\d _\a )=\g _\a $. By (7) we conclude $\forces _\LL \,
``\dot f :\h (n)^V \rightarrow \c ^V$ is onto''. 

A similar argument works for Miller forcing. $\qed $  

\sm

Combining Theorem 2.2 with Con$(\h (n+1)<\h (n))$ from $\S 1$ we
obtain the following:

\sm

{\bolds Corollary 2.3.} {\it For every $n<\o $, it is consistent that
both of Laver and Miller forcing collapse the continuum (strictly)
below $\h (n)$.}

\vfill \eject

{\gross References}
\Smallskip

\item{[Ba]}{J.E.Baumgartner, Iterated forcing, in: Surveys in set
theory, A.R.D. Mathias (ed.), London Math. Soc. Lect. Notes Ser. 8,
Cambridge Univ. Press, Cambridge (1983), 1--59}

\item{[BaPeSi]}{B. Balcar, J. Pelant, P.Simon, The space of
ultrafilters on $N$ covered by nowhere dense sets, Fund. Math. 110
(1980), 11--24}

\item{[Go]}{M. Goldstern, Tools for your forcing construction, in:
Israel Math. Conf. Proc. 6, H. Judah (ed.) (1993), 305--360} 

\item{[GoJoSp]}{M. Goldstern, M. Johnson and O. Spinas, 
Towers on trees, Proc. AMS 122 (1994),
557--564.} 

\item{[GoReShSp]}{M. Goldstern, M. Repick\`y, S. Shelah and O. Spinas,
On tree ideals, Proc. AMS 123 (1995), 1573--1581.}

\item{[JuSh]}{H. Judah and S. Shelah, Souslin forcing, J. Symb. Logic
53/4 (1988), 1188--1207.}

\item{[Mt]}{A.R.D. Mathias, Happy families, Ann. Math. Logic 12
(1977), 59--111.}

\item{[Shb]}{S. Shelah, Proper forcing, Lecture Notes in Math., vol.
942, Springer}

\item{[ShSp]}{S. Shelah and O. Spinas, The
distributivity number of $\P (\omega
)$/fin and its square, Trans. AMS, to appear.}

\bye